\documentclass[preprint,12pt]{elsarticle}
\usepackage{amssymb}
\usepackage{amsthm}
\usepackage{amsmath}
\usepackage{epic}
\newtheorem{thm}{Theorem}[section]
\newtheorem{cor}[thm]{Corollary}
\newtheorem{lem}[thm]{Lemma}
\newtheorem{pro}[thm]{Proposition}

\newtheorem{definition}[thm]{Definition}

\journal{~}

\begin{document}

\begin{frontmatter}

\title{Spectral properties of general hypergraphs}
\author[label1]{Changjiang Bu}\ead{buchangjiang@hrbeu.edu.cn}
\author[label1]{Jiang Zhou}
\author[label2]{Lizhu Sun}

\address[label1]{College of Science, Harbin Engineering University, Harbin 150001, PR China}
\address[label2]{School of Science, Harbin Institute of Technology, Harbin 150001, PR China}

\begin{abstract}
In this paper, we investigate spectral properties of the adjacency tensor, Laplacian tensor and signless Laplacian tensor of general hypergraphs (including uniform and non-uniform hypergraphs). We obtain some bounds for the spectral radius of general hypergraphs in terms of vertex degrees, and give spectral characterizations of odd-bipartite hypergraphs.
\end{abstract}


\end{frontmatter}

\section{Introduction}
Let $V(H)$ and $E(H)$ denote the vertex set and the edge set of a hypergraph $H$, respectively. A hypergraph $G$ satisfying $V(G)\subseteq V(H),E(G)\subseteq E(H)$ is called a \textit{sub-hypergraph} of $H$. If $G$ is a sub-hypergraph of $H$ and $G\neq H$, then $G$ is said to be a \textit{proper sub-hypergraph} of $H$. The \textit{degree} of a vertex $i$ of $H$ is defined as $d_i=|E_i|$, where $E_i$ denotes the set of edges containing $i$. The \textit{rank} and \textit{co-rank} of $H$ is the maximum and minimum cardinality of an edge in $H$, respectively \cite{Bretto}. Let $r(H)$ and $cr(H)$ denote the rank and co-rank of $H$, respectively. If $r(H)=cr(H)=k$, then $H$ is called $k$-\textit{uniform}. $2$-uniform hypergraphs are ordinary graphs. A \textit{path of length $l$} in a hypergraph $H$ is defined to be an alternating sequence $u_1e_1u_2\cdots u_le_lu_{l+1}$, where $u_1,\ldots,u_{l+1}$ are distinct vertices of $H$, $e_1,\ldots,e_l$ are distinct edges of $H$ and $u_i,u_{i+1}\in e_i$ for $i=1,\ldots,l$. If there exists a path between any two vertices of $G$, then $G$ is called \textit{connected}.

Hypergraphs and tensors are generalizations of graphs and matrices, and there is a natural one-to-one correspondence between uniform hypergraphs and tensors. In 2005, the concept of eigenvalues of tensors was posed by Qi \cite{Qi05} and Lim \cite{Lim}, independently. Due to the developments of spectral theory of tensors, there have been many attempts to extend spectral graph theory to hypergraphs via tensors. We first introduce some concepts and notations on tensor eigenvalues.

An order $m$ dimension $n$ tensor $\mathcal{A}=(a_{i_1i_2\cdots i_m})$ is a multidimensional array with $n^m$ entries ($i_j\in\{1,\ldots,n\},j=1,\ldots,m$). $\mathcal{A}$ is called \textit{symmetric} if $a_{i_1i_2\cdots i_k}=a_{i_{\sigma(1)}i_{\sigma(2)}\cdots i_{\sigma(k)}}$ for any permutation $\sigma$ on $\{1,\ldots,k\}$. For $\mathcal{A}=(a_{i_1i_2\cdots i_m})\in\mathbb{C}^{n\times n\times\cdots\times n}$ and $x=\left({x_1 ,\ldots ,x_n}\right)^\mathrm{T}\in\mathbb{C}^n,$ $\mathcal{A}x^{m-1}$ is a vector in $\mathbb{C}^n$ whose $i$-th component is
\begin{align*}
(\mathcal{A}x^{m-1})_i=\sum\limits_{i_2,\ldots,i_m=1}^na_{ii_2\cdots i_m}x_{i_2}\cdots x_{i_m}.
\end{align*}
A number $\lambda\in\mathbb{C}$ is called an \textit{eigenvalue} of $\mathcal{A}$, if there exists a nonzero vector $x\in\mathbb{C}^n$ such that $\mathcal{A}x^{m-1}=\lambda x^{[m-1]}$, where $x^{\left[ {m - 1} \right]}  = \left( {x_1^{m - 1},\ldots,x_n^{m - 1} } \right)^\mathrm{T}$. In this case, $x$ is an \textit{eigenvector} of $\mathcal{A}$ associated with $\lambda$. If $\lambda$ is a real eigenvalue with a real eigenvector, then $\lambda$ is called an \textit{H-eigenvalue} of $\mathcal{A}$. An H-eigenvalue with positive eigenvector is called an \textit{H$^{++}$-eigenvalue}. Let $\lambda(\mathcal{A})$ denote the largest H-eigenvalue of $\mathcal{A}$. The \textit{spectral radius} of $\mathcal{A}$ is defined as $\rho(\mathcal{A})=\max\{|\lambda|:\lambda\in\sigma(\mathcal{A})\}$, where $\sigma(\mathcal{A})$ is the set of all eigenvalues of $\mathcal{A}$.

In 2012, Cooper and Dutle \cite{Cooper12} investigated spectral properties of uniform hypergraphs via adjacency tensor. The \textit{adjacency tensor} of a $k$-uniform hypergraph $H$ with $n$ vertices, denoted by $\mathcal{A}_H$, is an order $k$ dimension $n$ symmetric tensor with entries
\begin{eqnarray*}
a_{i_1i_2\cdots i_k}=\begin{cases}\frac{1}{(k-1)!}~~~~~~~\{i_1,\ldots,i_k\}\in E(G),\\
0~~~~~~~~~~~~~\mbox{otherwise}.\end{cases}
\end{eqnarray*}
In 2014, Qi \cite{Qi14} defined the \textit{Laplacian tensor} and \textit{signless Laplacian tensor} of $H$ as $\mathcal{L}_H=\mathcal{D}_H-\mathcal{A}_H$ and $\mathcal{Q}_H=\mathcal{D}_H-\mathcal{A}_H$, respectively, where $\mathcal{D}_H$ is the diagonal tensor of vertex degrees of $H$. Recently, the research on spectral properties of $\mathcal{A}_H$, $\mathcal{L}_H$ and $\mathcal{Q}_H$ has attracted extensive attention [3,4,6,9,12-16].

In many literatures on spectral hypergraph theory, only uniform hypergraphs are considered. There are few work on spectral properties of non-uniform hypergraphs. In order to investigate spectra of non-uniform hypergraphs, we first extend the concept of the adjacency tensor, Laplacian tensor and signless Laplacian tensor from uniform hypergraphs to general hypergraphs as follows.
\begin{definition}\label{def1.1}
The \textit{adjacency tensor} of a hypergraph $H$ with $r(H)=k$, denoted by $\mathcal{A}_H$, is an order $k$ dimension $|V(H)|$ symmetric tensor with entries
\begin{eqnarray*}
a_{i_1i_2\cdots i_k}=\begin{cases}\frac{1}{(k-1)!}~~~~~\{i_1,\ldots,i_k\}\in E(H),\\
\frac{(k-s+1)!}{k!}~~\{i_1,\ldots,i_k\}=\{j_1^{(k-s+1)},j_2,\ldots,j_s\},~\{j_1,\ldots,j_s\}\in E(H),\\
0~~~~~~~~~~~\mbox{otherwise},\end{cases}
\end{eqnarray*}
where $s<k$ and $j_1^{(k-s+1)}$ means that the multiplicity of $j_1$ is $k-s+1$. Let $\mathcal{D}_H$ denote the order $k$ dimension $|V(H)|$ diagonal tensor whose diagonal entries are vertex degrees of $H$. The tensors $\mathcal{L}_H=\mathcal{D}_H-\mathcal{A}_H$ and $\mathcal{Q}_H=\mathcal{D}_H+\mathcal{A}_H$ are the \textit{Laplacian tensor} and the \textit{signless Laplacian tensor} of $H$, respectively.
\end{definition}

The following is an example for the tensor representations of a non-uniform hypergraph.

\vspace{3mm}
\noindent
\textbf{Example.} Let $H$ be a hypergraph whose vertex set and edge set are $V(H)=\{1,2,3,4,5\}$ and $E(H)=\{1234,45\}$, respectively. Then $r(H)=4$ and $cr(H)=2$. By Definition \ref{def1.1}, $\mathcal{A}_H=(a_{i_1i_2i_3i_4})$ is a symmetric tensor of dimension $5$, where $a_{i_1i_2i_3i_4}=\frac{1}{3!}=\frac{1}{6}$ if $\{i_1,i_2,i_3,i_4\}=\{1,2,3,4\}$, $a_{i_1i_2i_3i_4}=\frac{3!}{4!}=\frac{1}{4}$ if $\{i_1,i_2,i_3,i_4\}\in\{\{4,4,4,5\},\{4,5,5,5\}\}$, and $a_{i_1i_2i_3i_4}=0$ otherwise.

\vspace{3mm}
In this paper, we investigate spectral properties of the adjacency tensor, Laplacian tensor and signless Laplacian tensor of general hypergraphs. We obtain some bounds for the spectral radius of general hypergraphs in terms of vertex degrees, and give spectral characterizations of odd-bipartite hypergraphs.
\section{Preliminaries}

For a tensor $\mathcal{A}=(a_{i_1\cdots i_m })\in\mathbb{C}^{n\times\cdots\times n}$, we associate with $\mathcal{A}$ a digraph $\Gamma_{\mathcal{A}}$ as follows. The vertex set of $\Gamma_{\mathcal{A}}$ is $V(\mathcal{A})=\left\{1,\ldots,n\right\}$, the arc set of $\Gamma_{\mathcal{A}}$ is $E(\mathcal{A }) = \{(i,j)|a_{ii_2\cdots i_m}\neq0,j\in\{i_2,\ldots,i_m\}\neq\{i,\ldots,i\}\}$. $\mathcal{A}$ is said to be \textit{weakly irreducible} if $\Gamma_{\mathcal{A}}$ is strongly connected \cite{Bu,Pearson14}. We can get the following lemma from Definition \ref{def1.1}.
\begin{lem}\label{pro3.2}
For any hypergraph $H$, the following are equivalent:\\
(1) $H$ is connected.\\
(2) $\mathcal{A}_H$ is weakly irreducible.\\
(3) $\mathcal{L}_H$ is weakly irreducible.\\
(4) $\mathcal{Q}_H$ is weakly irreducible.
\end{lem}
A tensor $\mathcal{A}$ is said to be \textit{nonnegative} if all entries of $\mathcal{A}$ are nonnegative.
\begin{lem}\label{lem2.0}\textup{\cite{Qi14}}
Let $\mathcal{A}$ be a weakly irreducible nonnegative tensors of dimension $n$. Then $\rho(\mathcal{A})=\lambda(\mathcal{A})>0$ is the unique $H^{++}$-eigenvalue of $\mathcal{A}$.
\end{lem}

\begin{lem}\label{lem2.-1}\textup{\cite{LiShaoQi}}
Let $\mathcal{A}$ be a symmetric nonnegative tensors of order $k$ and dimension $n$. Then
\begin{eqnarray*}
\rho(\mathcal{A})=\lambda(\mathcal{A})=\max\left\{x^\textsc{T}(\mathcal{A}x^{k-1})|x\in\mathbb{R}^n_+,\sum_{i=1}^nx_i^k=1\right\}.
\end{eqnarray*}
Furthermore, $x\in\mathbb{R}^n_+$ with $\sum_{i=1}^nx_i^k=1$ is an eigenvector of $\mathcal{A}$ corresponding to $\rho(\mathcal{A})=\lambda(\mathcal{A})$ if and only if $\rho(\mathcal{A})=\lambda(\mathcal{A})=x^\textsc{T}(\mathcal{A}x^{k-1})$.
\end{lem}

For two nonnegative tensors $\mathcal{A}=(a_{i_1i_2\ldots i_m})$ and $\mathcal{B}=(b_{i_1i_2\ldots i_m})$ of dimension $n$, $\mathcal{B}\leqslant\mathcal{A}$ means that $b_{i_1i_2\cdots i_m}\leqslant a_{i_1i_2\cdots i_m}$ for all $i_1,i_2,\ldots,i_m\in\{1,\ldots,n\}$.
\begin{lem}\label{lem2.1}\textup{\cite{Kannan}}
Let $\mathcal{A}$ and $\mathcal{B}$ be two nonnegative tensors of order $m$ and dimension $n$, and $\mathcal{A}$ is weakly irreducible. If $\mathcal{B}\leqslant\mathcal{A}$ and $\mathcal{B}\neq\mathcal{A}$, then $\rho(\mathcal{B})<\rho(\mathcal{A})$.
\end{lem}

For $\mathcal{A}=(a_{i_1\cdots i_m})\in\mathbb{C}^{n\times\cdots\times n}$, the \textit{principal subtensor} of $\mathcal{A}$ with respect to $\alpha\subseteq\{1,\ldots,n\}$ is defined as $\mathcal{A}[\alpha]=(a_{i_1\cdots i_m})\in\mathbb{C}^{|\alpha|\times\cdots\times|\alpha|}$, where $i_1,\ldots,i_m\in\alpha$. If $|\alpha|<n$, then $\mathcal{A}[\alpha]$ is called \textit{proper principal subtensor} of $\mathcal{A}$.
\begin{lem}\label{lem2.5}\textup{\cite{Kannan}}
Let $\mathcal{A}$ be a weakly irreducible nonnegative tensors of dimension $n$. If $\mathcal{A}[\alpha]$ is a proper principal subtensor of $\mathcal{A}$, then $\rho(\mathcal{A}[\alpha])<\rho(\mathcal{A})$.
\end{lem}
In \cite{Shao 2013}, Shao defined the following product of tensors, which is a generalization of the matrix multiplication.
\begin{definition}\label{def2.1}\textup{\cite{Shao 2013}}
Let $\mathcal{A}$ and $\mathcal{B}$ be order $m\geqslant2$ and order $k\geqslant1$, dimension $n$ tensors, respectively. The product $\mathcal{A}\mathcal{B}$ is the following tensor $\mathcal{C}$ of order $(m-1)(k-1)+1$ and dimension $n$ with entries
\begin{eqnarray*}
c_{i\alpha_1\ldots \alpha_{m-1}}=\sum_{i_2,\ldots,i_m\in[n]}a_{ii_2\ldots i_m}b_{i_2\alpha_1}\cdots b_{i_m\alpha_{m-1}},
\end{eqnarray*}
where $i\in\{1,\ldots,n\},\alpha_1,\ldots,\alpha_{m-1}\in\{1,\ldots,n\}^{k-1}$.
\end{definition}

\begin{lem}\label{lem2.7}\textup{\cite{Zhou}}
Let $\mathcal{A}=(a_{i_1\cdots i_k})$ be an order $k\geqslant2$ dimension $n$ tensor, and let $P=(p_{ij}),Q=(q_{ij})$ be $n\times n$ matrices. Then
\begin{eqnarray*}
(P\mathcal{A}Q)_{i_1\cdots i_k}=\sum_{j_1,\ldots,j_k\in[n]}a_{j_1\cdots j_k}p_{i_1j_1}q_{j_2i_2}\cdots q_{j_ki_k}.
\end{eqnarray*}
\end{lem}
By Theorems 2.3 and 2.5 in \cite{Shao 2013}, we have the following lemma.
\begin{lem}\label{lem2.8}
Let $\mathcal{A}$ and $\mathcal{B}$ be two order $m$ dimension $n$ tensors. If there exists nonsingular diagonal real matrix $D$ such that $\mathcal{B}=D^{-(m-1)}\mathcal{A}D$, then $\mathcal{A}$ and $\mathcal{B}$ have the same spectrum and H-spectrum.
\end{lem}

\section{Main results}
Let $H$ be a hypergraph with $n$ vertices and $r(H)=k$. For $x\in\mathbb{C}^n$ and a vertex subset $\alpha$ of $H$, let $x^\alpha=\prod_{i\in\alpha}x_i$. By Definition \ref{def1.1}, we have
\begin{eqnarray*}
(\mathcal{A}_Hx^{k-1})_i&=&\sum_{i_2,\ldots,i_k=1}^na_{ii_2\cdots i_k}x_{i_2}\cdots x_{i_k}\\
&=&\frac{1}{k}\sum_{e\in E_i}[(k-|e|)x^{e\setminus\{i\}}x_i^{k-|e|}+x^{e\setminus\{i\}}\sum_{j\in e}x_j^{k-|e|}],~i=1,\ldots,n,\\
x_i(\mathcal{A}_Hx^{k-1})_i&=&\frac{1}{k}\sum_{e\in E_i}[(k-|e|)x^ex_i^{k-|e|}+x^e\sum_{j\in e}x_j^{k-|e|}],~i=1,\ldots,n,\\
\mathcal{A}_Hx^k&=&x^\textsc{T}(\mathcal{A}x^{k-1})=\sum_{i_1,\ldots,i_k=1}^na_{i_1i_2\cdots i_k}x_{i_1}\cdots x_{i_k}=\sum_{e\in E(H)}x^e\sum_{j\in e}x_j^{k-|e|}.
\end{eqnarray*}
Moreover, if $H$ is $k$-uniform, then
\begin{eqnarray*}
x^\textsc{T}(\mathcal{A}x^{k-1})=k\sum_{e\in E(H)}x^e,~(\mathcal{A}_Hx^{k-1})_i=\sum_{e\in E_i}x^{e\setminus\{i\}},~x_i(\mathcal{A}_Hx^{k-1})_i=\sum_{e\in E_i}x^e
\end{eqnarray*}
for $i=1,\ldots,n$.

From the above equalities, we obtain the following result.
\begin{pro}\label{pro3.1}
Let $H$ be a hypergraph with $n$ vertices and $r(H)=k$. If $\lambda$ is an eigenvalue of $H$ with an eigenvector $x$, then
\begin{eqnarray*}
\lambda x_i^{k-1}&=&\frac{1}{k}\sum_{e\in E_i}[(k-|e|)x^{e\setminus\{i\}}x_i^{k-|e|}+x^{e\setminus\{i\}}\sum_{j\in e}x_j^{k-|e|}],~i=1,\ldots,n,\\
\lambda x_i^k&=&\frac{1}{k}\sum_{e\in E_i}[(k-|e|)x^ex_i^{k-|e|}+x^e\sum_{j\in e}x_j^{k-|e|}],~i=1,\ldots,n,\\
\lambda\sum_{i=1}^nx_i^k&=&\sum_{e\in E(H)}x^e\sum_{j\in e}x_j^{k-|e|}.
\end{eqnarray*}
Moreover, if $H$ is $k$-uniform, then
\begin{eqnarray*}
\lambda\sum_{j=1}^nx_j^k=k\sum_{e\in E(H)}x^e,~\lambda x_i^{k-1}=\sum_{e\in E_i}x^{e\setminus\{i\}},~\lambda x_i^k=\sum_{e\in E_i}x^e,~i=1,\ldots,n.
\end{eqnarray*}
\end{pro}

\begin{thm}
Let $H$ be a connected hypergraph. Then the following hold:\\
(1) $\rho(\mathcal{A}_H)=\lambda(\mathcal{A}_H)>0$ is the unique $H^{++}$-eigenvalue of $\mathcal{A}_H$.\\
(2) $\rho(\mathcal{Q}_H)=\lambda(\mathcal{Q}_H)>0$ is the unique $H^{++}$-eigenvalue of $\mathcal{Q}_H$.\\
(3) If $G$ is a proper sub-hypergraph of $H$ with $r(G)=r(H)$, then
\begin{eqnarray*}
\rho(\mathcal{A}_G)<\rho(\mathcal{A}_H),~\rho(\mathcal{Q}_G)<\rho(\mathcal{Q}_H).
\end{eqnarray*}
\end{thm}
\begin{proof}
Since $H$ is connected, by Lemma \ref{pro3.2}, $\mathcal{A}_H$ and $\mathcal{Q}_H$ are weakly irreducible. By Lemma \ref{lem2.0}, we know that (1) and (2) hold.

Since $G$ is a proper sub-hypergraph of $H$ with $r(G)=r(H)$, one of the following cases holds:

(a) $\mathcal{A}_G\leqslant\mathcal{A}_H,\mathcal{Q}_G\leqslant\mathcal{Q}_H$ and $\mathcal{A}_G\neq\mathcal{A}_H,\mathcal{Q}_G\neq\mathcal{Q}_H$.

(b) $\mathcal{A}_G$ and $\mathcal{Q}_G$ are proper principal subtensors of $\mathcal{A}_H$ and $\mathcal{Q}_H$, respectively.

Since $\mathcal{A}_H$ is weakly irreducible, by Lemmas \ref{lem2.1} and \ref{lem2.5}, we know that (3) holds.
\end{proof}


\begin{thm}
Let $H$ be a hypergraph with average degree $\overline{d}$. Then $\rho(\mathcal{A}_H)\geqslant\overline{d}$, with equality if and only if $H$ is regular.
\end{thm}
\begin{proof}
Let $x=(\frac{1}{\sqrt[k]{n}},\ldots,\frac{1}{\sqrt[k]{n}})^\textsc{T}$, where $k=r(H),n=|V(H)|$. By Lemma \ref{lem2.-1}, we have
\begin{eqnarray*}
\rho(\mathcal{A}_H)\geqslant x^\textsc{T}(\mathcal{A}_Hx^{k-1})=\sum_{e\in E(H)}x^e\sum_{j\in e}x_j^{k-|e|}=\frac{1}{n}\sum_{e\in E(H)}|e|=\overline{d},
\end{eqnarray*}
with equality if and only if $\mathcal{A}_Hx^{k-1}=\rho(\mathcal{A}_H)x^{[k-1]}$, i.e., $H$ is regular.
\end{proof}

\begin{thm}
Let $H$ be a connected hypergraph with maximum degree $\Delta$. Then $\rho(\mathcal{A}_H)\leqslant\Delta$, with equality if and only if $H$ is regular.
\end{thm}
\begin{proof}
Let $x=(x_1,\ldots,x_n)^\textsc{T}$ be the positive eigenvector corresponding to $\rho(H)$, and let $x_i=\max_{1\leqslant j\leqslant n}x_j$. By Proposition \ref{pro3.1}, we have
\begin{eqnarray*}
\rho(\mathcal{A}_H)x_i^k&=&\frac{1}{k}\sum_{e\in E_i}[(k-|e|)x^ex_i^{k-|e|}+x^e\sum_{j\in e}x_j^{k-|e|}]\leqslant d_ix_i^k,
\end{eqnarray*}
with equality if and only if $x_j=x_i$ for any $j\in e\in E_i$. Since $H$ is connected, we have $\rho(\mathcal{A}_H)\leqslant\Delta$, with equality if and only if $H$ is regular.
\end{proof}

\begin{thm}
Let $H$ be a connected hypergraph with $r(H)=k$. Then
\begin{eqnarray*}
\rho(\mathcal{A}_H)\leqslant\max\left\{\sqrt[k]{d_{i_1}^{k-s+1}d_{i_2}\cdots d_{i_s}}~|~\{i_1,\ldots,i_s\}\in E(H),d_{i_1}\geqslant\cdots\geqslant d_{i_s}\right\}.
\end{eqnarray*}
\end{thm}
\begin{proof}
Let $x$ be the positive eigenvector corresponding to $\rho(\mathcal{A}_H)$. Suppose that $x_{i_1}\cdots x_{i_k}=\max_{(\mathcal{A}_H)_{j_1\cdots j_k}\neq0}x_{j_1}\cdots x_{j_k}$. By Proposition \ref{pro3.1}, we have
\begin{eqnarray*}
\rho(\mathcal{A}_H)x_{i_j}^k=\frac{1}{k}\sum_{e\in E_{i_j}}[(k-|e|)x^ex_{i_j}^{k-|e|}+x^e\sum_{l\in e}x_l^{k-|e|}]\leqslant d_{i_j}x_{i_1}\cdots x_{i_k}
\end{eqnarray*}
for all $j=1,\ldots,k$. Then
\begin{eqnarray*}
\rho(\mathcal{A}_H)^k\prod_{j=1}^kx_{i_j}^k&\leqslant&d_{i_1}\cdots d_{i_k}(x_{i_1}\cdots x_{i_k})^k,\\
\rho(\mathcal{A}_H)&\leqslant&\sqrt[k]{d_{i_1}\cdots d_{i_k}}.
\end{eqnarray*}
By Definition \ref{def1.1}, we have
\begin{eqnarray*}
\rho(\mathcal{A}_H)\leqslant\max\left\{\sqrt[k]{d_{i_1}^{k-s+1}d_{i_2}\cdots d_{i_s}}~|~\{i_1,\ldots,i_s\}\in E(H),d_{i_1}\geqslant\cdots\geqslant d_{i_s}\right\}.
\end{eqnarray*}
\end{proof}
We can obtain the following result from the above Theorem.
\begin{cor}
Let $H$ be a connected $k$-uniform hypergraph. Then
\begin{eqnarray*}
\rho(\mathcal{A}_H)\leqslant\max_{\{i_1,\ldots,i_k\}\in E(H)}\sqrt[k]{d_{i_1}d_{i_2}\cdots d_{i_k}}.
\end{eqnarray*}
\end{cor}

\noindent
\textbf{Remark.} In \cite{Yuan}, Yuan et al. proved that
\begin{eqnarray*}
\rho(\mathcal{A}_H)\leqslant\max_{e\in E(H)}\max_{\{i,j\}\subseteq e}\sqrt{d_id_j}.
\end{eqnarray*}
Corollary 3.6 is a refinement of this upper bound.

\vspace{3mm}
A hypergraph $H$ is called \textit{odd-bipartite}, if its vertex set has a partition $V(H)=V_1\cup V_2$ such that each edge of $H$ contains odd number of vertices in $V_1$ and odd number of vertices in $V_2$. Spectral characterizations of odd-bipartite uniform hypergraphs are given in \cite{HuQi-DAM,ShaoShanWu}. We extend these work to general odd-bipartite hypergraphs as follows.
\begin{thm}\label{thm3.5}
Let $G$ be a connected hypergraph. Then the following are equivalent:\\
(1) $G$ is odd-bipartite.\\
(2) $\mathcal{A}_G$ and $-\mathcal{A}_G$ have the same spectrum and H-spectrum.\\
(3) $-\rho(\mathcal{A}_G)$ is an H-eigenvalue of $\mathcal{A}_G$.
\end{thm}
\begin{proof}
(1)$\Rightarrow$(2). If $G$ is odd-bipartite, then by Lemma \ref{lem2.7}, there exists a diagonal matrix $P$ with diagonal entries $\pm1$ such that $\mathcal{A}_G=-P^{-(k-1)}\mathcal{A}_GP$, where $k=r(G)$. By Lemma \ref{lem2.8}, we know that $\mathcal{A}_G$ and $-\mathcal{A}_G$ have the same spectrum and H-spectrum.

(2)$\Rightarrow$(3). If $\mathcal{A}_G$ and $-\mathcal{A}_G$ have the same H-spectrum, then $-\rho(\mathcal{A}_G)$ is an H-eigenvalue of $\mathcal{A}_G$.

(3)$\Rightarrow$(1). By Lemma 1.2 in \cite{ShaoShanWu}, there exits real diagonal matrix $P=diag(x_1,\ldots,x_n)$ such that $\mathcal{A}_G=-P^{-(k-1)}\mathcal{A}_GP$ and $|x_1|=\cdots|x_n|=1$. By Lemma \ref{lem2.7}, we have
\begin{eqnarray*}
(\mathcal{A}_G)_{i_1i_2\cdots i_k}=-(\mathcal{A}_G)_{i_1i_2\cdots i_k}x_{i_1}^{-(k-1)}x_{i_2}\cdots x_{i_k}.
\end{eqnarray*}
For any $(\mathcal{A}_G)_{i_1i_2\cdots i_k}\neq0$, we have
\begin{eqnarray*}
x_{i_1}x_{i_2}\cdots x_{i_k}=-x_{i_1}^k=\cdots=-x_{i_k}^k.
\end{eqnarray*}
If $k$ is odd, then $x_{i_1}=\cdots=x_{i_k}$, a contradiction to $x_{i_1}x_{i_2}\cdots x_{i_k}=-x_{i_1}^k$. So $k$ is even and $x_{i_1}x_{i_2}\cdots x_{i_k}=-1$ for any $(\mathcal{A}_G)_{i_1i_2\cdots i_k}\neq0$. Let $V_1=\{u|u\in V(G),x_u=-1\}$. By Definition 1.1, we know that each edge $e$ of $G$ contains odd number of vertices in $V_1$ and $|e|$ is even. Hence $G$ is odd-bipartite.
\end{proof}

\begin{thm}
Let $G$ be a connected hypergraph. Then the following are equivalent:\\
(1) $G$ is odd-bipartite.\\
(2) $\mathcal{L}_G$ and $\mathcal{Q}_G$ have the same spectrum and H-spectrum.\\
(3) $0$ is an H-eigenvalue of $\mathcal{Q}_G$.
\end{thm}
\begin{proof}
(1)$\Rightarrow$(2). If $G$ is odd-bipartite, then by Lemma \ref{lem2.7}, there exists a diagonal matrix $P$ with diagonal entries $\pm1$ such that $\mathcal{L}_G=P^{-(k-1)}\mathcal{Q}_GP$, where $k=r(G)$. By Lemma \ref{lem2.8}, we know that $\mathcal{L}_G$ and $\mathcal{Q}_G$ have the same spectrum and H-spectrum.

(2)$\Rightarrow$(3). Let $x=(1,\ldots,1)^\textsc{T}$. Since $\mathcal{L}_Gx^{k-1}=0$, $0$ is always an H-eigenvalue of $\mathcal{L}_G$. If $\mathcal{L}_G$ and $\mathcal{Q}_G$ have the same H-spectrum, then $0$ is an H-eigenvalue of $\mathcal{Q}_G$.

(3)$\Rightarrow$(1). Let $x=(x_1,\ldots,x_n)^\textsc{T}$ be a real eigenvector of $\mathcal{Q}_G$ corresponding to the H-eigenvalue $0$. Then
\begin{eqnarray*}
(\mathcal{Q}_Gx^{k-1})_i=d_ix_i^{k-1}+(\mathcal{A}_Gx^{k-1})_i=0,~i=1,\ldots,n,\\
\frac{1}{k}\sum_{e\in E_i}[(k-|e|)x^{e\setminus\{i\}}x_i^{k-|e|}+x^{e\setminus\{i\}}\sum_{j\in e}x_j^{k-|e|}]=-d_ix_i^{k-1},~i=1,\ldots,n.
\end{eqnarray*}
let $|x_j|=\max_{1\leqslant i\leqslant n}|x_i|$. Then
\begin{eqnarray*}
-d_ix_j^k=\frac{1}{k}\sum_{e\in E_j}[(k-|e|)x^ex_j^{k-|e|}+x^e\sum_{i\in e}x_i^{k-|e|}],\\
d_i|x_j^k|\leqslant\frac{1}{k}\sum_{e\in E_j}[(k-|e|)|x_j^k|+|e||x_j^k|]=d_i|x_j^k|.
\end{eqnarray*}
Hence $x^ex_i^{k-|e|}=-x_i^k$ for any $i\in e$ and any $e\in E(G)$. If $k$ is odd, then $x_1=\cdots=x_n$, a contradiction to $x^ex_i^{k-|e|}=-x_i^k$. So $k$ is even and $x^ex_i^{k-|e|}=-|x_j|^k$ for any $i\in e$ and any $e\in E(G)$. Let $V_1=\{u|u\in V(G),x_u=-|x_j|\}$, then each edge $e$ of $G$ contains odd number of vertices in $V_1$ and $|e|$ is even. Hence $G$ is odd-bipartite.
\end{proof}

\noindent
\textbf{Acknowledgements.}

\vspace{3mm}
This work is supported by the National Natural Science Foundation of China (No. 11371109), the Natural Science Foundation of the Heilongjiang Province (No. QC2014C001) and the Fundamental Research Funds for the Central Universities.

\end{document}